\documentclass[11pt]{amsart}
\usepackage[margin=1.2in]{geometry}
\usepackage{amssymb,amsmath,amsthm,mathtools}
\usepackage{enumitem}
\usepackage{array}
\usepackage[colorlinks,linkcolor=blue,citecolor=blue,urlcolor=blue]{hyperref}

\newtheorem{theorem}{Theorem}[section]
\newtheorem{proposition}[theorem]{Proposition}
\newtheorem{lemma}[theorem]{Lemma}
\newtheorem{corollary}[theorem]{Corollary}
\newtheorem{fact}[theorem]{Fact}
\theoremstyle{definition}
\newtheorem{definition}[theorem]{Definition}
\theoremstyle{remark}
\newtheorem{remark}[theorem]{Remark}
\newtheorem{question}[theorem]{Question}

\newcommand{\N}{\mathbb{N}}
\newcommand{\F}{\mathbb{F}}
\newcommand{\cant}{2^{\N}}
\newcommand{\str}{2^{<\N}}
\newcommand{\cA}{\mathcal{A}}
\newcommand{\supp}{\operatorname{supp}}
\newcommand{\rank}{\operatorname{rank}}
\newcommand{\jump}[1]{0^{(#1)}}
\newcommand{\MLR}{\mathsf{MLR}}
\newcommand{\KLR}{\mathsf{KLR}}
\newcommand{\CR}{\mathsf{CR}}
\newcommand{\SR}{\mathsf{SR}}
\newcommand{\WR}{\mathsf{WR}}
\newcommand{\Low}{\operatorname{Low}}
\newcommand{\E}{\mathbb{E}}
\numberwithin{equation}{section}

\title[Base dependence and lowness for KL-randomness]{Base dependence and lowness for\\ Kolmogorov--Loveland randomness}

\author{Bj{\o}rn Kjos-Hanssen}

\date{27 November 2026}

\begin{document}

\begin{abstract}
We derive three consequences of the recent separation of Kolmogorov--Loveland
(KL) randomness from Martin-L\"of (ML) randomness by Liu and Meng. First,
KL-randomness of a real number depends on the base in which the real is
written: there is a real whose base-$4$ expansion is KL-random while its
binary expansion is not. The converse implication, from base $2$ to base $4$,
holds. The proof combines the separation with the splitting theorem of
Merkle, Miller, Nies, Reimann and Stephan. We run the Liu--Meng construction
against the betting strategies that read the two bits of each base-$4$ digit
together, and place its parity constraints alternately in the even half and
in the odd half of the binary expansion. Second, we bound the complexity of
the separation itself: there is a sequence that is KL-random and not
ML-random and is computable in the fifth jump $\jump5$, and a measure
concentrated on such sequences whose cylinder probabilities are rational and
computable in $\jump4$. Third, the only sets that are low for KL-randomness
are the computable ones. This settles a question left open by Nies, who
showed that such sets are $K$-trivial. The new ingredient is that they are
hyperimmune-free: if $A$ has hyperimmune degree, then some KL-random sequence
is not computably random relative to $A$, and in fact not Schnorr random
relative to $A$. It also follows that every KL-random sequence is Schnorr
random relative to $A$ if and only if $A$ is computably traceable. The first
two results, and the hyperimmune-freeness of the sets that are low for
KL-randomness, have been formalized in Lean~4.
\end{abstract}

\maketitle

\section{Introduction}

ML-randomness, computable randomness and Schnorr randomness of a real number
do not depend on the base used to write it; for computable randomness this is
a theorem of Brattka, Miller and Nies \cite{BMN}. Whether the same holds for
KL-randomness is a question of A.~Nies, recorded as open in
\cite[Section 9]{Rute}. As long as KL-randomness might have coincided with
ML-randomness \cite{MN}, a positive answer was possible for a trivial reason.
Liu and Meng \cite{LM} have now separated the two notions: they construct a
probability measure $\rho$ on $\cant$, supported on sequences failing one
fixed ML-test, such that $\rho$-almost every sequence is KL-random. Their
theorem has been verified in Lean \cite{klr}. After this separation, any
failure of base invariance has to be witnessed by a KL-random real that is
not ML-random. We show that such a failure occurs already for the bases $2$
and $4$.

For $X\in\cant$ let $W_X\in4^{\N}$ be given by $W_X(n)=2X(2n)+X(2n+1)$, so
that $X$ and $W_X$ are the binary and base-$4$ expansions of the same real
$0.X$. Write $X=X_0\oplus X_1$, where $X_c(n)=X(2n+c)$.

\begin{theorem}\label{thm:main}
There is a real that is KL-random in base $4$ and not KL-random in base $2$.
More precisely, there is a Borel probability measure $\rho'$ on $\cant$ such
that
\begin{enumerate}[label=\textup{(\roman*)}]
\item for every $X\in\supp\rho'$, neither $X_0$ nor $X_1$ is ML-random, and
hence $X$ is not KL-random in base $2$;
\item for $\rho'$-almost every $X$, the sequence $W_X$ is KL-random in base
$4$.
\end{enumerate}
\end{theorem}

\begin{proposition}\label{prop:easy}
Let $b\ge 2$ and $k\ge 1$. If a real is KL-random in base $b$, then it is
KL-random in base $b^k$.
\end{proposition}

\begin{corollary}\label{cor:strict}
The set of reals that are KL-random in base $2$ is a proper subset of the set
of reals that are KL-random in base $4$.
\end{corollary}

Our second result concerns the separation of \cite{LM} itself. The
construction there makes several choices that are not effective: a list of
all total computable strategies, a subsequence along which certain functions
converge weakly, and the blocks of parity constraints. In particular the
measure $\rho$ cannot be computable (Remark~\ref{rem:notcomp}). We show that
all choices can be made within the arithmetical hierarchy. Write $\jump{n}$
for the $n$th jump of the empty set.

\begin{theorem}\label{thm:eff}
\begin{enumerate}[label=\textup{(\alph*)}]
\item There is a Borel probability measure $\rho$ on $\cant$ such that
$\rho$-almost every $X$ is KL-random and not ML-random, the cylinder
probabilities $\rho[\sigma]$ are rational, and the function
$\sigma\mapsto\rho[\sigma]$ is computable in $\jump4$.
\item There is a sequence $X\le_T\jump5$ that is KL-random and not
ML-random.
\end{enumerate}
\end{theorem}

We do not know whether these bounds are optimal; see
Questions~\ref{q:four} and~\ref{q:measure}.

Our third result concerns lowness. Let
$\MLR\subseteq\KLR\subseteq\CR\subseteq\SR$ be the classes of ML-random,
KL-random, computably random and Schnorr random sequences in $\cant$. If a
class $\mathcal D$ is given by a relativizable definition, let $\mathcal D^A$
be its relativization to an oracle $A$. Following \cite{KNS,N05}, for
$\mathcal C\subseteq\mathcal D$ let
\[
\Low(\mathcal C,\mathcal D)=\{A:\mathcal C\subseteq\mathcal D^A\},\qquad
\Low(\mathcal C)=\Low(\mathcal C,\mathcal C).
\]
A stronger oracle yields a smaller class $\mathcal D^A$, so membership in
$\Low(\mathcal C,\mathcal D)$ expresses computational weakness and is closed
downward under Turing reducibility. For
$\mathcal C,\mathcal D\in\{\MLR,\CR,\SR\}$ these classes were determined by
Nies \cite{N05} and by Kjos-Hanssen, Nies and Stephan \cite{KNS}, building on
\cite{TZ,BN}; see also \cite{BM}. They are the entries of
Table~\ref{tab:low} outside the $\KLR$ row and column.

The sets in $\Low(\KLR)$ are the sets that are \emph{low for KL-randomness}.
Nies showed that every such set is low for $K$, hence $K$-trivial
\cite[Corollary 5.4]{N05}, and noted that it was unknown whether a
non-computable one exists \cite[p.~277]{N05}. There is none.

\begin{theorem}\label{thm:lowklr}
A set is low for KL-randomness if and only if it is computable. In fact
$\Low(\KLR)=\Low(\KLR,\CR)$ is the class of computable sets.
\end{theorem}

A function $g\colon\N\to\N$ is \emph{dominated} by $h$ if $g(n)\le h(n)$ for
all but finitely many $n$. A set has \emph{hyperimmune degree} if it computes
a function that is not dominated by any computable function, and is
\emph{hyperimmune-free} otherwise \cite{MM}.

The proof of Theorem~\ref{thm:lowklr} has the form of Nies's proof that only
the computable sets are low for computable randomness
\cite[Theorem 5.7]{N05}. A set in $\Low(\KLR,\CR)$ is in $\Low(\MLR,\CR)$ and
is therefore low for $K$ \cite[Theorem 5.2]{N05}. So it is $K$-trivial, hence
$\Delta^0_2$ \cite{Ch}, and a hyperimmune-free $\Delta^0_2$ set is computable
\cite{MM}. What remains is to show that every set in $\Low(\KLR,\CR)$ is
hyperimmune-free. For $\Low(\CR)$ this is a theorem of Bedregal and Nies
\cite{BN}. For KL-randomness it follows from the next theorem, which is the
new ingredient, and this consequence has been formalized in Lean.

\begin{theorem}\label{thm:hyp}
Let $g\colon\N\to\N$ be a function that is not dominated by any computable
function. There are a Borel probability measure $\rho$ on $\cant$ and a
Schnorr test relative to $g$ such that $\rho$-almost every sequence is
KL-random and every sequence in $\supp\rho$ fails the test. In particular, if
$A$ has hyperimmune degree, then some KL-random sequence is not Schnorr
random relative to $A$, and hence not computably random relative to $A$.
\end{theorem}

By the last clause, every set in $\Low(\KLR,\CR)$ is hyperimmune-free. The
clause about Schnorr randomness determines one more class.

\begin{corollary}\label{cor:klrsr}
The following are equivalent for a set $A$.
\begin{enumerate}[label=\textup{(\roman*)}]
\item Every KL-random sequence is Schnorr random relative to $A$.
\item $A$ is computably traceable.
\end{enumerate}
Thus $\Low(\KLR,\SR)=\Low(\CR,\SR)=\Low(\SR)$.
\end{corollary}

One further entry of Table~\ref{tab:low} does not need the separation and is
implicit in \cite[Section 5]{N05}.

\begin{proposition}\label{prop:mlrklr}
$\Low(\MLR,\KLR)$ is the class of $K$-trivial sets.
\end{proposition}

\begin{table}[ht]
\renewcommand{\arraystretch}{1.25}
\begin{tabular}{@{}lllll@{}}
\hline
$\mathcal C\ \backslash\ \mathcal D$ & $\MLR$ & $\KLR$ & $\CR$ & $\SR$\\
\hline
$\MLR$ & $K$-trivial & \textbf{$\boldsymbol{K}$-trivial} & $K$-trivial & c.e.\ traceable\\
$\KLR$ & & \textbf{computable} & \textbf{computable} & \textbf{computably traceable}\\
$\CR$ & & & computable & computably traceable\\
$\SR$ & & & & computably traceable\\
\hline
\end{tabular}
\medskip
\caption{The classes $\Low(\mathcal C,\mathcal D)$. The entries in bold are
Theorem~\ref{thm:lowklr}, Corollary~\ref{cor:klrsr} and
Proposition~\ref{prop:mlrklr}.}\label{tab:low}
\end{table}

In Table~\ref{tab:low} the $\KLR$ row agrees with the $\CR$ row, and the
$\KLR$ column agrees with the $\CR$ column.

Theorem~\ref{thm:hyp} strengthens the separation of \cite{LM}. Let $A$ be
c.e.\ traceable and of hyperimmune degree, for instance a non-computable
$K$-trivial set. The theorem yields a KL-random sequence $Z$ that is not
Schnorr random relative to $A$, and $Z$ is not ML-random because
$A\in\Low(\MLR,\SR)$. Conversely, if $\MLR=\KLR$ held, then $\Low(\KLR,\SR)$
would equal $\Low(\MLR,\SR)$, the class of c.e.\ traceable sets. So
Corollary~\ref{cor:klrsr} could not have been proved without separating the
two notions.

\subsection*{Idea of the proof of Theorem~\ref{thm:main}}
The $n$th base-$4$ digit of $0.X$ is the pair $(X_0(n),X_1(n))$. A base-$4$
strategy is the same thing as a binary strategy that always reads the bits
$2n$ and $2n+1$ one right after the other (Lemma~\ref{lem:sim}). By the
splitting theorem \cite[Theorem 12]{MMNRS}, a KL-random binary sequence
cannot have two non-ML-random halves. The proof of that theorem uses
strategies that scan one half far ahead before betting on the other, which is
exactly what a base-$4$ strategy cannot do. We therefore run the construction
of \cite{LM} against the pair-respecting strategies only, and put the parity
constraints of even stages on even positions and those of odd stages on odd
positions. Both halves then fail an ML-test. The proof of \cite{LM} that
almost every sequence defeats the listed strategies still applies, because a
pair-respecting scan rule queries position $2n$ if and only if it queries
$2n+1$, so one family of Borel partitions serves both halves.

\subsection*{Idea of the proof of Theorem~\ref{thm:hyp}}
The measure of \cite{LM} is concentrated on sequences that satisfy infinitely
many systems of parity equations. The $j$th system has rank $r_j$ and is the
output of a program $p_j$ of length much less than $r_j$, and the ML-test
enumerates the solution sets of all high-rank systems with short programs.
We arrange in addition that $p_j$ halts within $g(r_j)$ steps. The part of
the test that waits only $g(r)$ steps for systems of rank $r$ is a Schnorr
test relative to $g$, and it still covers the support of the measure. Two
properties of the construction make this possible. The rank $r_j$ may be
chosen among all sufficiently large numbers. And after a small change to the
construction (Lemma~\ref{lem:horizon}) the running time of $p_j$ is bounded
by a computable function $f_j$ of the rank. Since $g$ is not dominated by
$f_j$, some admissible rank $r$ satisfies $g(r)\ge f_j(r)$. This is the
pattern of the argument of Bedregal and Nies that the sets in
$\Low(\CR,\SR)$ are hyperimmune-free \cite{BN}, \cite[Lemma 4.9]{KNS}.

\subsection*{Formalization}
Theorems~\ref{thm:main} and~\ref{thm:eff} and Corollary~\ref{cor:strict} have
been formalized in Lean~4 with Mathlib \cite{klr}. So has the part of
Theorem~\ref{thm:hyp} on which Theorem~\ref{thm:lowklr} rests: if $A$ has
hyperimmune degree, then some KL-random sequence is not computably random
relative to $A$, and hence every set that is low for KL-randomness is
hyperimmune-free. Section~\ref{sec:lean} lists the correspondence between the
statements below and the formal ones, and says what is not formalized. For
Theorem~\ref{thm:eff} we give an outline of the proof in
Section~\ref{sec:eff} and refer to the formalization for the details. The
proof of Theorem~\ref{thm:hyp} is given in full in Section~\ref{sec:low}.

\section{Betting in base $b$ and pair-respecting strategies}\label{sec:bet}

Let $b\ge2$ and write $b^{<\N}$ for the finite strings over $\{0,\dots,b-1\}$.
A \emph{base-$b$ scan rule} is a total computable $\pi\colon b^{<\N}\to\N$
that chooses distinct positions along every chain of prefixes. A
\emph{base-$b$ strategy} is a base-$b$ scan rule together with a uniformly
computable capital function $d\colon b^{<\N}\to\mathbb{R}_{\ge0}$ with
\[
d(\epsilon)=1,\qquad b\,d(v)=\sum_{c<b}d(vc).
\]
On input $W\in b^{\N}$ the transcript is $v_0=\epsilon$,
$v_{t+1}=v_t\,W(\pi(v_t))$, and the strategy \emph{succeeds} on $W$ if
$\sup_t d(v_t)=\infty$. A sequence $W\in b^{\N}$ is \emph{KL-random in base
$b$} if no base-$b$ strategy succeeds on it, and a real in $[0,1]$ is
KL-random in base $b$ if its base-$b$ expansion is. For $b=2$ this is the
definition used in \cite{LM} and \cite{MMNRS}, in terms of total strategies.

\begin{definition}
A binary scan rule $\pi$ is \emph{pair-respecting} if for every
$v\in\str$ of even length, $\pi(v)$ is even and $\pi(v0)=\pi(v1)=\pi(v)+1$.
A binary strategy is pair-respecting if its scan rule is.
\end{definition}

Along every run, a pair-respecting rule queries $2a_0,2a_0+1,2a_1,2a_1+1,\dots$
for distinct $a_0,a_1,\dots$. In particular, for every even $a$ and every
infinite run,
\begin{equation}\label{eq:pair}\tag{P}
\text{$\pi$ queries $a$}\iff\text{$\pi$ queries $a+1$.}
\end{equation}

\begin{lemma}\label{lem:sim}
For $X\in\cant$ the following are equivalent:
\begin{enumerate}[label=\textup{(\roman*)}]
\item some base-$4$ strategy succeeds on $W_X$;
\item some pair-respecting binary strategy succeeds on $X$.
\end{enumerate}
\end{lemma}

\begin{proof}
Identify a binary string $v$ of even length with the base-$4$ string
$\bar v$ obtained by pairing its bits.

(i)$\Rightarrow$(ii). Let $(\pi_4,d_4)$ be a base-$4$ strategy. For $v$ of
even length put $\pi(v)=2\pi_4(\bar v)$, $\pi(vb)=2\pi_4(\bar v)+1$,
\[
d(v)=d_4(\bar v),\qquad d(vb)=\tfrac12\bigl(d_4(\bar v\,(2b))+d_4(\bar v\,(2b+1))\bigr).
\]
Then $d(v0)+d(v1)=\frac12\sum_{c<4}d_4(\bar vc)=2d(v)$ and
$d(vb0)+d(vb1)=2d(vb)$, so $(\pi,d)$ is a pair-respecting binary strategy.
Its capital at even times along $X$ is the capital of $(\pi_4,d_4)$ along
$W_X$.

(ii)$\Rightarrow$(i). Let $(\pi,d)$ be pair-respecting. Put
$\pi_4(\bar v)=\pi(v)/2$ and $d_4(\bar v)=d(v)$. Fairness holds because
$\sum_{b,b'}d(vbb')=2\sum_bd(vb)=4d(v)$. Since $d(vb)\le2d(v)$, the capital
of $(\pi,d)$ is unbounded along $X$ only if it is unbounded at even times.
\end{proof}

\begin{proof}[Proof of Proposition~\ref{prop:easy}]
For $b=2$, $k=2$ this is (i)$\Rightarrow$(ii) of Lemma~\ref{lem:sim}: a
pair-respecting strategy is a binary strategy. In general a base-$b^k$
strategy is simulated in base $b$ by reading the $k$ base-$b$ digits of a
block consecutively and splitting the bet on the block into $k$ conditional
bets, exactly as above.
\end{proof}

\begin{remark}[partial strategies]\label{rem:partial}
All strategies in this paper are total. In base $2$, partial computable
strategies define the same notion of KL-randomness \cite{Merkle}; see
\cite{MMNRS}. We do not treat partial strategies in base $4$ here.
\end{remark}

\section{The splitting theorem}\label{sec:split}

\begin{theorem}[{\cite[Theorem 12]{MMNRS}}]\label{thm:split}
Let $Z\subseteq\N$ be computable, infinite and co-infinite. If $A\in\cant$ is
KL-random, then at least one of the two subsequences of $A$ selected by $Z$
and by $\N\setminus Z$ is ML-random.
\end{theorem}

\begin{corollary}\label{cor:split}
If neither $X_0$ nor $X_1$ is ML-random, then $X=X_0\oplus X_1$ is not
KL-random in base $2$.
\end{corollary}

The strategies behind Theorem~\ref{thm:split} use non-monotonicity only
between the two halves: one half is scanned, without betting, far enough to
extract a bound that is then used to bet on the other half. Such a strategy
separates the bits $2n$ and $2n+1$ by arbitrarily long stretches of the run,
so it is not pair-respecting, and Lemma~\ref{lem:sim} yields no base-$4$
strategy from it. The analogue of Theorem~\ref{thm:split} for base-$4$
sequences concerns splittings of the \emph{digit} positions and says nothing
about the two bit-coordinates of the digits.

\section{The Liu--Meng stage construction with confined pools}\label{sec:LM}

We use the notation of \cite{LM}. $\lambda$ is fair-coin measure on $\cant$
and $\kappa$ is fair-coin measure on the space $\cant$ of answer sequences.
For a finite $S\subset\N$, a binary matrix $B$ with columns indexed by $S$
and full row rank $r$, and $y\in\F_2^r$,
\[
E(S,B,y)=\{X\in\cant: B(X{\restriction}S)=y\},\qquad\lambda(E(S,B,y))=2^{-r}.
\]
Programs are plain binary programs for a fixed universal interpreter, as in
\cite[Section 3]{LM}.

Fix total computable binary scan rules $(\pi_i)_{i\ge1}$ and Borel sets
$H_i\subseteq\cant$ of answer sequences, with complements $L_i$.

\begin{definition}\label{def:good}
Let $\ell\ge1$, $\eta>0$ and $0<\tau<1$. A finite set $G\subset\N$ is
\emph{$(\ell,\eta,\tau)$-good} if for every $i\le\ell$,
\begin{align}
\kappa\{Z\in H_i:\ \pi_i\text{ misses a position of }G\}&<\eta,\tag{4}\\
\kappa\{Z\in L_i:\ \pi_i\text{ misses fewer than }\tau|G|\text{ positions of }G\}&<\eta.\tag{5}
\end{align}
For $\cA\subseteq\N$ we say that $(H_i)$ \emph{supplies pools in $\cA$} if for
all $\ell\ge1$, rational $\eta>0$ and $z\in\N$ there are a rational
$0<\tau<1$ and a program that, for every sufficiently large $M$, outputs an
$(\ell,\eta,\tau)$-good set $G\subseteq\cA\cap[z,\infty)$ with $|G|=M$.
\end{definition}

The equation numbers (4) and (5) are those of \cite{LM}. In this language,
\cite[Lemma 3]{LM} says that for every sequence of total computable scan
rules there is $(H_i)$ supplying pools in $\N$.

\begin{proposition}[stage construction with confined pools]\label{prop:stage}
Let $(\sigma_i)_{i\ge1}$ be any sequence of total computable binary
strategies, normalized as in \cite[Lemma 2]{LM}, with scan rules $\pi_i$.
Let $(H_i)$ be Borel sets and $\cA_1,\cA_2,\dots\subseteq\N$ be such that
$(H_i)$ supplies pools in $\cA_j$ for every $j$. Then there are pairwise
disjoint finite sets $S_j\subseteq\cA_j$, full-row-rank matrices $B_j$ with
$r_j$ rows and columns indexed by $S_j$, syndromes $y_j$, and programs $p_j$
outputting $(S_j,B_j,y_j)$, such that
\begin{enumerate}[label=\textup{(\alph*)}]
\item $r_j\ge2|p_j|+j+1$ for all $j\ge1$;
\item if $\rho$ is the law under which $X{\restriction}S_j$ is a uniform
solution of $B_jx=y_j$, independently over $j$, and all other bits are
independent and fair, then $\supp\rho=\bigcap_jE(S_j,B_j,y_j)$ and, for every
$i$, the capital of $\sigma_i$ is bounded $\rho$-almost surely.
\end{enumerate}
\end{proposition}

\begin{proof}
This is the proof of \cite[Sections 5 and 6]{LM}, read with two
observations.

\emph{The list of strategies is arbitrary.} The construction of \cite{LM}
starts from a list of all total computable strategies, but uses the list only
through its first $j$ entries at stage $j$. That the list is exhaustive
enters once, at the very end, to pass from ``every listed strategy has
bounded capital almost surely'' to KL-randomness.

\emph{The pool lemma is used only through its conclusion.} At stage $j$,
\cite[Section 5]{LM} applies Lemma 3 with $\ell=j$,
$\eta=e_j/(100j2^K)$ and a position bound $z$. From then on it uses the pool
program, the fraction $\tau$, the bounds (4) and (5) for the pool $G$ of size
$M=r^8$, and the fact that the sets $H_i$ are the same at all stages. The
candidate supports are subsets of $G$, so $S_j\subseteq\cA_j$ when the pool
program outputs subsets of $\cA_j$. The estimates (6)--(9) of \cite{LM}, the
choice of horizon and syndrome, the description bound (a), and the whole of
\cite[Section 6]{LM} make no reference to where the pool lies.
\end{proof}

\section{Pools for pair-respecting scan rules}\label{sec:pools}

\begin{lemma}[even pools]\label{lem:even}
For every sequence $(\pi_i)_{i\ge1}$ of total computable scan rules there are
Borel sets $H_i$ such that $(H_i)$ supplies pools in $2\N$.
\end{lemma}

\begin{proof}
Let $\varepsilon$ be the scan rule that reads the even positions
$0,2,4,\dots$ in increasing order, whatever the answers. Put
$\pi'_1=\pi'_2=\varepsilon$ and $\pi'_{i+2}=\pi_i$, and run the proof of
\cite[Lemma 3]{LM} for the family $(\pi'_i)$. With its weights $w_i=2^{-i}$,
the weighted query probability $p(a)$ satisfies
\[
p(a)\ge w_1+w_2=\tfrac34\ \text{ for even }a,\qquad
p(a)\le\sum_{i\ge3}w_i=\tfrac14\ \text{ for odd }a.
\]
Hence $p^*=\limsup_ap(a)\ge\frac34$, and the positions $a_1<a_2<\cdots$ with
$p(a_k)\to p^*$ chosen in that proof are even from some point on. The
indicator that $\varepsilon$ queries $a_k$ is then identically $1$, so its
limit is $g'_1=1$ almost everywhere and the set $H'_1=\{g'_1=1\}$ has
$\kappa(H'_1)=1$.

Let $(H'_i)$ be the sets produced by that proof, and put $H_i=H'_{i+2}$.
Given $\ell$, $\eta$ and $z$, apply the conclusion of \cite[Lemma 3]{LM} for
$(\pi'_i)$ with $\ell+2$, $\min(\eta,\frac12)$ and $z$, obtaining $\tau$ and a
program. Let $G$ be one of its outputs. If $G$ contained an odd position,
$\varepsilon$ would miss it along every answer sequence, and the event in (4)
for $i=1$ would be all of $H'_1$, of measure $1>\frac12$. So $G\subseteq2\N$.
The bounds (4) and (5) for the indices $3,\dots,\ell+2$ say that $G$ is
$(\ell,\eta,\tau)$-good for $(\pi_i)$ and $(H_i)$.
\end{proof}

The auxiliary rule $\varepsilon$ serves only to steer the pools. It is not
among the strategies that the final construction defeats.

\begin{lemma}[parity transfer]\label{lem:transfer}
Suppose every $\pi_i$ is pair-respecting. If $(H_i)$ supplies pools in $2\N$,
then the same $(H_i)$ supplies pools in $2\N+1$.
\end{lemma}

\begin{proof}
Given $\ell,\eta,z$, take $\tau$ and the program for $2\N$, and let the new
program output $G+1=\{a+1:a\in G\}$. By \eqref{eq:pair}, along every answer
sequence $\pi_i$ misses $a\in G$ if and only if it misses $a+1$. So the events
in (4) and (5) are the same for $G$ and for $G+1$.
\end{proof}

\section{Proof of Theorem~\ref{thm:main}}

\subsection*{Tests for the two halves}
For $c\in\{0,1\}$ and a system $(S,B,y)$ with $S\subseteq2\N+c$ let
$(S,B,y)^{(c)}$ be the system on $\{(a-c)/2:a\in S\}$ with the same matrix
and syndrome, so that
\begin{equation}\label{eq:half}\tag{H}
X\in E(S,B,y)\iff X_c\in E\bigl((S,B,y)^{(c)}\bigr).
\end{equation}
For $m\ge1$ let $V^c_m$ be the union of the sets $E((S,B,y)^{(c)})$ over all
programs $p$ that output a system $(S,B,y)$ with $S\subseteq2\N+c$ and
$\rank B\ge2|p|+m+1$.

\begin{lemma}\label{lem:tests}
For each $c\in\{0,1\}$, $(V^c_m)_{m\ge1}$ is an ML-test. If programs $p_j$
output systems $(S_j,B_j,y_j)$ with $S_j\subseteq2\N+(j\bmod2)$ and ranks
$r_j\ge2|p_j|+j+1$ for all $j\ge1$, and $X\in\bigcap_jE(S_j,B_j,y_j)$, then
$X_0\in\bigcap_mV^0_m$ and $X_1\in\bigcap_mV^1_m$.
\end{lemma}

\begin{proof}
The sets $V^c_m$ are uniformly effectively open, as for the test of
\cite[Section 3]{LM}. A system of rank $r$ has
$\lambda(E((S,B,y)^{(c)}))=2^{-r}$, and each program outputs at most one
system, so
\[
\lambda(V^c_m)\le\sum_{d\ge0}2^d\,2^{-2d-m-1}=2^{-m}.
\]
For the second statement fix $c$ and $m$, and take $j\ge m$ with
$j\equiv c\pmod2$. Then $r_j\ge2|p_j|+m+1$, so
$E((S_j,B_j,y_j)^{(c)})\subseteq V^c_m$, and $X_c$ lies in this set by
\eqref{eq:half}.
\end{proof}

The test uses the original programs $p_j$; no recoding of a system to its
half, and no change in program length, is involved.

\subsection*{The construction}
Let $(\sigma_i)_{i\ge1}$ list all pair-respecting total computable binary
strategies, each normalized by \cite[Lemma 2]{LM}. The normalized strategy
has the same scan rule, so it is again pair-respecting, and it succeeds
wherever the original does. Let $\pi_i$ be the scan rule of $\sigma_i$.

By Lemma~\ref{lem:even} there are Borel sets $H_i$ such that $(H_i)$ supplies
pools in $2\N$, and by Lemma~\ref{lem:transfer} also in $2\N+1$. Put
$\cA_j=2\N+(j\bmod2)$. Apply Proposition~\ref{prop:stage} and let $\rho'$ be
the resulting law, with blocks $(S_j,B_j,y_j)$, $S_j\subseteq\cA_j$, and
programs $p_j$.

\subsection*{Proof of (i)}
Let $X\in\supp\rho'=\bigcap_jE(S_j,B_j,y_j)$. By
Proposition~\ref{prop:stage}(a) and Lemma~\ref{lem:tests}, $X_0$ fails the
test $(V^0_m)$ and $X_1$ fails the test $(V^1_m)$. So neither half of $X$ is
ML-random, and by Corollary~\ref{cor:split}, $X$ is not KL-random in base
$2$.

\subsection*{Proof of (ii)}
By Proposition~\ref{prop:stage}(b), $\rho'$-almost surely every $\sigma_i$ has
bounded capital. The $\sigma_i$ are countably many and cover all
pair-respecting total computable strategies up to normalization. Hence for
$\rho'$-almost every $X$ no pair-respecting strategy succeeds on $X$, and by
Lemma~\ref{lem:sim} no base-$4$ strategy succeeds on $W_X$.

\smallskip
For $\rho'$-almost every $X$ the real $0.X$ is therefore KL-random in base
$4$ and not in base $2$. Such an $X$ is not eventually constant, so both
expansions of the real are unique. This proves Theorem~\ref{thm:main}, and
Corollary~\ref{cor:strict} follows with Proposition~\ref{prop:easy}. \qed

\begin{remark}\label{rem:notcomp}
The measure $\rho'$ is not computable. If it were, the computable martingale
$\sigma\mapsto2^{|\sigma|}\rho'[\sigma]$ would succeed on every sequence in
the support, and a left-to-right strategy is pair-respecting. The same
remark applies to the measure $\rho$ of \cite{LM}. More generally, no
computable measure $\mu$ is concentrated on the sequences that are KL-random
and not ML-random: such a $\mu$ is singular with respect to $\lambda$, so the
computable martingale $\sigma\mapsto2^{|\sigma|}\mu[\sigma]$ succeeds
$\mu$-almost surely. Theorem~\ref{thm:eff}(a) shows that such a measure can
nevertheless be computable in $\jump4$.
\end{remark}

\section{Why the strategies must respect pairs}\label{sec:check}

The construction cannot be carried out against the list of all binary
strategies: it would produce a KL-random sequence with two non-ML-random
halves, contradicting Theorem~\ref{thm:split}. It is instructive to see
where it stops.

Lemma~\ref{lem:even} holds for arbitrary scan rules, and so does its mirror
image for odd pools, obtained from an auxiliary rule reading the odd
positions. But the two lemmas then produce two families of limits
$g_i^{0},g_i^{1}$ and two partitions $(H_i^{0},L_i^{0})$ and
$(H_i^{1},L_i^{1})$, and Lemma~\ref{lem:transfer} fails: for the rule that
reads exactly the odd positions in increasing order, $g^{1}=1$ and $g^{0}=0$
everywhere. Proposition~\ref{prop:stage} needs one family $(H_i)$ for all
stages. In \cite[Section 6]{LM}, on the event that the scan has no rank drop
on block $j$ but has one on a later block $k$, transcripts in $H_i$ are
handled by (8) for stage $j$ and transcripts in $L_i$ by (9) for stage $k$.
With two partitions, a transcript in $L_i^{c(j)}\cap H_i^{c(k)}$ with
$c(j)\ne c(k)$ is covered by neither estimate. That is the transcript of a
strategy which reads all odd positions and no even ones: it passes through
all later odd blocks, meets rank drops there, and may afterwards return to
the even block $j$ with what it has learned. This is the shape of the
strategies behind Theorem~\ref{thm:split}.

Conversely, in the construction of \cite{LM} one subsequence $(a_k)$ serves
all scan rules. Let $Z$ be computable, infinite and co-infinite. For the two
rules that read exactly $Z$, respectively exactly $\N\setminus Z$,
convergence of $\int g_{i,a_k}\,d\kappa$ forces $(a_k)$ to lie eventually on
one side of $Z$, and then (4) for the rule reading that side forces every
pool from some stage on to lie on that side. So all but finitely many blocks
of \cite{LM} lie on one side of each computable splitting, in agreement with
Theorem~\ref{thm:split}. For pair-respecting rules the splitting into even and
odd bits is invisible to this mechanism, by \eqref{eq:pair}.

\section{The complexity of a separating sequence}\label{sec:eff}

We outline the proof of Theorem~\ref{thm:eff}. It is an effective version of
the construction of \cite{LM}, in which every non-effective choice is
replaced by a search relative to a finite jump. The complete proof is the
formalization \cite{klr}. Table~\ref{tab:jumps} records where the jumps are
spent.

\begin{table}[ht]
\renewcommand{\arraystretch}{1.25}
\begin{tabular}{@{}>{\raggedright\arraybackslash}p{0.66\textwidth}c@{}}
\hline
Object & Computable in\\
\hline
A uniform family of total strategies covering all strategies & no oracle\\
Normalized strategies; positions $a_k$ and the limits of the query coordinates & $\jump2$\\
$L^1$ names of the weak limits $g_i$ & $\jump3$\\
Clopen approximations in measure of the sets $H_i$ & $\jump4$\\
The sequence of blocks $(S_j,B_j,y_j)$ and programs $p_j$ & $\jump4$\\
The cylinder probabilities $\rho[\sigma]$ & $\jump4$\\
A KL-random sequence in $\supp\rho$ & $\jump5$\\
\hline
\end{tabular}
\medskip
\caption{Upper bounds in the proof of Theorem~\ref{thm:eff}.}\label{tab:jumps}
\end{table}

\subsection*{A uniform family of strategies}
The construction of \cite{LM} starts from a list of all total computable
strategies, which cannot be effective. Instead one uses a uniformly
computable sequence $(\tau_e)_{e\ge1}$ of total strategies with rational
capital, obtained by simulating all programs with bounded resources, two
simulators for each program, as in the passage from partial to total
strategies \cite{Merkle}. For every total computable strategy $\sigma$ there
are $e_0,e_1$ such that $\tau_{e_0}$ or $\tau_{e_1}$ succeeds wherever
$\sigma$ succeeds. Hence $X$ is KL-random if and only if no $\tau_e$ succeeds
on $X$. By Proposition~\ref{prop:stage}, whose proof works for an arbitrary
list, it suffices to run the construction against the normalizations
$\sigma_e$ of the $\tau_e$. In the formalization these are uniformly
computable in $\jump2$.

\subsection*{Limits of the query coordinates in $\jump2$}
The proof of \cite[Lemma 3]{LM} chooses positions $a_1<a_2<\cdots$ along which
the numbers $\int_{[v]}g_{i,a}\,d\kappa$ converge for all $i$ and all
cylinders $[v]$, and along which the weighted query probability $p(a)$ tends
to its limsup. Call these numbers, together with $p(a)$, the \emph{query
coordinates} of $a$. They are uniformly left-c.e.\ reals in $[0,1]$. A
\emph{reservoir} is a set of positions of the form $[z,\infty)$ intersected
with finitely many strict rational lower bounds on query coordinates.
Reservoirs are uniformly c.e.; whether a reservoir is infinite is a $\Pi_2$
question and whether it is empty is a $\Pi_1$ question. With the oracle
$\jump2$ one narrows, one coordinate at a time, a rational interval for the
limit of that coordinate: the lower bound is chosen so that the reservoir
stays infinite, and the upper bound is certified by the emptiness of a
reservoir with a higher threshold. This yields positions $a_k$, computable
in $\jump2$, along which every query coordinate converges, and rational
approximations to the limits with error $2^{-k}$, also computable in
$\jump2$.

The limit of $p(a_k)$ obtained in this way need not be the limsup of $p$.
The proof of \cite[Lemma 3]{LM} only needs the following substitute: for every
$s>0$ there is a reservoir that contains all sufficiently late $a_k$ and on
which $p$ is at most the limit plus $s$. The pool search is run inside such a
reservoir. A reservoir is determined by finitely many parameters and is
c.e., so the pool program is still an ordinary program, as the description
bound requires.

\subsection*{The sets $H_i$ in $\jump4$}
Let $g_i$ be the weak limit of $g_{i,a_k}$ and $H_i=\{g_i=1\}$. The integrals
of $g_i$ over cylinders are the limits named above, so they are computable in
$\jump2$. The averages of $g_i$ over the cylinders of length $d$ converge to
$g_i$ in $L^1(\kappa)$ as $d\to\infty$, with no rate known in advance. The
Cauchy condition for a sequence of rational step functions computable in an
oracle $A$ is $\Pi_1$ in $A$, so a modulus of convergence is found by a search
in $A'$. This gives $L^1$ names of the $g_i$ in $\jump3$. Applying the ramps
$x\mapsto\max(0,1-m(1-x))$, which converge to the indicator of $\{x=1\}$ on
$[0,1]$, and accelerating once more, one obtains clopen sets that approximate
$H_i$ in $\kappa$-measure to within $2^{-k}$, uniformly computable in
$\jump4$. This is not a procedure deciding membership in $H_i$.

\subsection*{The stages in $\jump4$}
At stage $j$ the construction must find a support, a matrix, a horizon and
a syndrome satisfying the estimates (6)--(9) of \cite{LM}, and a program
$p_j$ of length at most $(r_j-j-1)/2$ that outputs the block. Applying the
existence argument of \cite[Section 5]{LM} with half the allowed error shows
that there is a candidate meeting all estimates with a margin. Each estimate
is then a strict inequality between a rational number and a probability that
can be approximated, with a computable error bound, from the clopen
approximations of the sets $H_i$, from finite simulations of the scan rules,
from approximations, computable in $0'$, of the events that a rank drop ever
occurs, and from exact rational conditioning on the earlier blocks. So $\jump4$ can
search for a candidate together with a finite certificate, and the sequence
$j\mapsto(S_j,B_j,y_j,p_j)$ is computable in $\jump4$.

The oracle selects the programs $p_j$ but is not given to them: each $p_j$
is an ordinary program for the fixed interpreter. So every sequence in
$\bigcap_jE(S_j,B_j,y_j)$ fails the unrelativized ML-test of
\cite[Section 3]{LM}.

\subsection*{Proof of Theorem~\ref{thm:eff}(a)}
Let $\rho$ be the block law of the selected blocks. By
Proposition~\ref{prop:stage}, $\rho$-almost every $X$ defeats every
$\sigma_e$, hence every $\tau_e$, and is therefore KL-random; and every $X$
in the support fails an ML-test. The position bound of stage $j$ is chosen
at least $j$, so $\min S_j\ge j$. A cylinder $[\sigma]$ with $|\sigma|=n$
depends only on positions below $n$, so $\rho[\sigma]$ equals the
probability of $[\sigma]$ under the law of the first $n$ blocks. That
probability is a ratio of two counts of solutions of linear systems over
$\F_2$. Hence $\rho[\sigma]$ is rational and computable from the first $n$
blocks, and so in $\jump4$.

\subsection*{Proof of Theorem~\ref{thm:eff}(b)}
Let $A=\jump4$. For $e,b\ge1$ let $C(e,b)$ be the set of $X$ along which the
capital of $\tau_e$ never exceeds $b$. These sets are closed, uniformly
effectively so, they increase with $b$, and
$\rho(\bigcup_bC(e,b))=1$ for every $e$ by part (a). We build strings
$\sigma_0\subset\sigma_1\subset\cdots$ with $|\sigma_n|=n$ and bounds
$b_1,b_2,\dots$ such that
\[
\rho\Bigl([\sigma_n]\cap\bigcap_{e\le n}C(e,b_e)\Bigr)\ge4^{-n}.
\]
Given stage $n$, some $b_{n+1}$ keeps at least half of this measure inside
$C(n+1,b_{n+1})$, and then one of the two extensions of $\sigma_n$ by one bit
keeps at least half of what remains. The displayed condition says that every
clopen approximation from above of the closed set has $\rho$-measure at least
$4^{-n}$. These measures are rational and computable in $A$ by part (a), so
the condition is $\Pi_1$ in $A$, and $A'=\jump5$ can search for $b_{n+1}$ and
$\sigma_{n+1}$. Let $X=\bigcup_n\sigma_n$. Then $X\le_T\jump5$. Every
prefix of $X$ has positive measure inside each closed set $C(e,b_e)$, so $X$
lies in all of them and no $\tau_e$ succeeds on $X$: $X$ is KL-random. Every
prefix of $X$ has positive $\rho$-measure, so $X\in\supp\rho$ and $X$ is not
ML-random. \qed

\begin{remark}
The effective construction selects one particular instance of the
construction of \cite{LM}. It does not show that every measure obtained from
the choices allowed in \cite{LM} is computable in $\jump4$.
\end{remark}

\section{Lowness}\label{sec:low}

We prove Theorem~\ref{thm:hyp} and derive Theorem~\ref{thm:lowklr},
Corollary~\ref{cor:klrsr} and Proposition~\ref{prop:mlrklr}. The proof of
Theorem~\ref{thm:hyp} is one more variant of the construction of \cite{LM}.
Here the list of strategies is the full one and the pools are not confined.
Instead, the programs $p_j$ are made to halt quickly relative to $g$.

\subsection*{Background and derivations}
A \emph{Schnorr test relative to $A$} is a sequence
$(\mathcal O_m)_{m\ge1}$ of uniformly $\Sigma^0_1(A)$ classes with
$\lambda(\mathcal O_m)\le2^{-m}$ such that $\lambda(\mathcal O_m)$ is
computable relative to $A$, uniformly in $m$. A sequence \emph{fails} the
test if it lies in $\bigcap_m\mathcal O_m$, and is Schnorr random relative to
$A$ if it fails no such test. Requiring $\lambda(\mathcal O_m)=2^{-m}$
instead gives the same notion; see, e.g., \cite[Chapter 3]{N09}.

A set $A$ is \emph{c.e.\ traceable} if there is a computable function $u$
such that every function $\varphi\le_TA$ has a trace: a uniformly c.e.\
sequence $(F_n)$ of finite sets with $|F_n|\le u(n)$ and $\varphi(n)\in F_n$
for all $n$. It is \emph{computably traceable} if moreover the trace can be
chosen to be a computable sequence of canonical indices of finite sets.

\begin{fact}\label{fact:known}
\leavevmode
\begin{enumerate}[label=\textup{(\alph*)}]
\item\label{f:nies} $\Low(\MLR)=\Low(\MLR,\CR)$ is the class of sets that are
low for $K$ \cite[Theorem 5.2]{N05}, which is the class of $K$-trivial sets
\cite[Theorem 6.2]{N05}.
\item\label{f:ce} $\Low(\MLR,\SR)$ is the class of c.e.\ traceable sets
\cite[Theorem 3.1(I)]{KNS}.
\item\label{f:ct} $\Low(\CR,\SR)=\Low(\SR)$ is the class of computably
traceable sets \cite{TZ}, \cite[Theorem 3.1(II)]{KNS}.
\item\label{f:hif} A hyperimmune-free c.e.\ traceable set is computably
traceable \cite[Proposition 4.14]{KNS}.
\item\label{f:delta} $K$-trivial sets are $\Delta^0_2$ \cite{Ch}, and a
hyperimmune-free $\Delta^0_2$ set is computable \cite{MM}.
\end{enumerate}
\end{fact}

For every oracle $A$ we have
$\MLR^A\subseteq\KLR^A\subseteq\CR^A\subseteq\SR^A$. For the two inclusions
involving $\KLR^A$ it does not matter whether $\KLR^A$ is defined with
partial or with total $A$-computable strategies, and these inclusions are the
only properties of $\KLR^A$ used below.

\begin{proof}[Proof of Theorem~\ref{thm:lowklr} from Theorem~\ref{thm:hyp}]
Computable sets are in $\Low(\KLR)$. Since $\KLR^A\subseteq\CR^A$, every set
in $\Low(\KLR)$ is in $\Low(\KLR,\CR)$. Let $A\in\Low(\KLR,\CR)$, so that
every KL-random sequence is computably random relative to $A$. By
Theorem~\ref{thm:hyp}, $A$ is hyperimmune-free. Since $\MLR\subseteq\KLR$,
$A\in\Low(\MLR,\CR)$, so $A$ is $K$-trivial by
Fact~\ref{fact:known}\ref{f:nies}, hence $\Delta^0_2$. By
Fact~\ref{fact:known}\ref{f:delta}, $A$ is computable.
\end{proof}

This is the proof of \cite[Theorem 5.7]{N05}, with the theorem of Bedregal
and Nies \cite{BN} replaced by Theorem~\ref{thm:hyp}. The other steps are
Fact~\ref{fact:known}\ref{f:nies} and~\ref{f:delta}, which are classical. The
step taken from Theorem~\ref{thm:hyp} is formalized; see
Section~\ref{sec:lean}.

\begin{proof}[Proof of Corollary~\ref{cor:klrsr} from Theorem~\ref{thm:hyp}]
If $A$ is computably traceable then $\KLR\subseteq\SR\subseteq\SR^A$ by
Fact~\ref{fact:known}\ref{f:ct}. Conversely, let $A\in\Low(\KLR,\SR)$. Then
$A\in\Low(\MLR,\SR)$ because $\MLR\subseteq\KLR$, so $A$ is c.e.\ traceable by
Fact~\ref{fact:known}\ref{f:ce}. By Theorem~\ref{thm:hyp}, $A$ is
hyperimmune-free. Now apply Fact~\ref{fact:known}\ref{f:hif}.
\end{proof}

\begin{proof}[Proof of Proposition~\ref{prop:mlrklr}]
Since $\MLR^A\subseteq\KLR^A\subseteq\CR^A$,
\[
\Low(\MLR)\subseteq\Low(\MLR,\KLR)\subseteq\Low(\MLR,\CR),
\]
and the two outer classes coincide with the $K$-trivial sets by
Fact~\ref{fact:known}\ref{f:nies}.
\end{proof}

\subsection*{The stage estimates of Liu and Meng}
We use the notation of Section~\ref{sec:LM}. Let $(\sigma_i)_{i\ge1}$ list all
total computable binary strategies, each normalized by \cite[Lemma 2]{LM}.
Let $\pi_i$ be the scan rule of $\sigma_i$ and $h_i$ its signed wager
fraction, so that $\sum_th_i(v_t)^2\le\frac12$ along every run. Let $H_i$,
with complements $L_i$, be the Borel sets of answer sequences that
\cite[Lemma 3]{LM} provides for $(\pi_i)$. If $P$ is a probability measure on
inputs and $F$ is a Borel set of answer sequences, then $P(F)$ denotes the
$P$-probability that the answer sequence of the scan rule under consideration
lies in $F$. Thus $\lambda(F)=\kappa(F)$.

A \emph{system} is a triple $(S,B,y)$ as in Section~\ref{sec:LM}, and its
\emph{rank} is the number of rows of $B$. While a scan rule reads an input,
the \emph{remaining rank} of the system is the rank of the columns of $B$ at
the positions of $S$ that have not been queried. A query that lowers the
remaining rank is a \emph{rank drop}. Under $\lambda$ conditioned on
$E(S,B,y)$, the answer at a rank drop is determined by the earlier answers,
and every other answer is a fair bit given the earlier answers
\cite[Section 3]{LM}. Let $U_m$ be the union of all $E(S,B,y)$ such that some
program $p$ outputs the system $(S,B,y)$ and its rank is at least
$2|p|+m+1$. Then $(U_m)_{m\ge1}$ is the ML-test of \cite[Section 3]{LM}.

Suppose we are given systems $(S_j,B_j,y_j)$ of rank $r_j$, for $j\ge1$,
with pairwise disjoint supports, and numbers $V_j\in\N$, the \emph{horizons}.
Write
\[
E_j=E(S_j,B_j,y_j),\qquad P_j=\lambda(\,\cdot\mid E_1\cap\dots\cap E_j),
\qquad P_0=\lambda,\qquad e_j=2^{-j}/100,
\]
and let $\rho$ be the law of Proposition~\ref{prop:stage}(b). For the rule
$\pi_i$ and the block $j$:
\begin{itemize}
\item $R_{ij}(t)$ is the remaining rank of $(S_j,B_j,y_j)$ after $t$ queries;
\item $T_{ij}\in\{1,2,\dots\}\cup\{\infty\}$ is the index of the first query
that lowers $R_{ij}$;
\item $W_{ij}(t)$ is the sum of the wager fractions $|h_i|$ at the rank drops
on block $j$ among the first $t$ queries;
\item for a finite set $S$, $A_i(S)$ is the sum of $|h_i(v_t)|$ over all
queries with $\pi_i(v_t)\in S$. Thus $W_{ij}(t)\le A_i(S_j)$.
\end{itemize}
The estimates of \cite[Section 5]{LM} at stage $j$ are the following, for
every $i\le j$. As for (4) and (5), the numbers are those of \cite{LM}.
\begin{align}
\E_{P_j}W_{ij}(V_j)+P_j\bigl(0<R_{ij}(V_j)<r_j\bigr)&<4e_j, \tag{6}\\
P_j(V_j<T_{ij}<\infty)&\le2e_j, \tag{7}\\
P_{j-1}\bigl(H_i\cap\{T_{ij}=\infty\}\bigr)&<e_j, \tag{8}\\
P_a\bigl(L_i\cap\{T_{ij}<\infty\}\bigr)&<e_j\qquad(0\le a<j). \tag{9}
\end{align}

\begin{theorem}[{\cite[Section 6]{LM}}]\label{thm:LM6}
Suppose that for every $j\ge1$
\begin{enumerate}[label=\textup{(\alph*)}]
\item\label{avoid} no position that one of $\pi_1,\dots,\pi_j$ can query
among its first $V_j$ queries, along any answer sequence, belongs to a
support $S_k$ with $k>j$, and
\item \textup{(6)--(9)} hold for every $i\le j$.
\end{enumerate}
Then for every $i$ the capital of $\sigma_i$ is bounded $\rho$-almost surely.
Hence $\rho$-almost every sequence is KL-random.
\end{theorem}

Section 6 of \cite{LM} uses nothing about the blocks and the horizons beyond
these hypotheses. This is also how Proposition~\ref{prop:stage}(b) is
obtained.

We need the following form of \cite[Lemma 3]{LM}. It is what the proof given
there establishes.

\begin{lemma}[pools with horizons]\label{lem:poolV}
Fix $\ell\ge1$, a rational $0<\eta<1$ and $z\in\N$. There are a rational
$0<\tau<1$ and a program which, for every sufficiently large $M$, halts on
input $M$ and outputs a pair $(G,V)$ such that $G\subseteq[z,\infty)$ is
$(\ell,\eta,\tau)$-good with $|G|=M$, $V\in\N$, and for every $i\le\ell$
\begin{equation}\label{eq:high}\tag{4$^+$}
\kappa\{Z\in H_i:\ \pi_i\text{ does not query all of $G$ among its first $V$ queries}\}<\eta.
\end{equation}
\end{lemma}

\begin{proof}
The program in the proof of \cite[Lemma 3]{LM} enumerates pairs $(G,V)$ and
outputs $G$ from the first pair that passes three tests. We let it output the
pair. Every set $G$ output in this way is $(\ell,\eta,\tau)$-good, as shown
there. The first test requires
$\kappa\{Z\in H_i':v_{i,V}(Z)<1\}<3\delta$ for $i\le\ell$, where $v_{i,V}(Z)$
is the fraction of $G$ that $\pi_i$ reads among its first $V$ queries, $H_i'$
is a clopen set with $\kappa(H_i\triangle H_i')<\delta$, and $4\delta<\eta$.
Hence $\kappa\{Z\in H_i:v_{i,V}(Z)<1\}<\delta+3\delta<\eta$, which
is~\eqref{eq:high}.
\end{proof}

In \cite{LM}, (4) is derived from \eqref{eq:high}: a position that is never
queried has not been queried among the first $V$ queries.

\subsection*{Candidates at stage \texorpdfstring{$j$}{j}}
Suppose the stages $a<j$ have been carried out. Put $K=\sum_{a<j}r_a$. Then
$P_a\le2^K\lambda$ for all $a<j$, hence $P_a(F)\le2^K\kappa(F)$ for every
Borel set $F$ of answer sequences of a fixed scan rule. Choose $z$ above all
earlier supports and above every position that one of $\pi_1,\dots,\pi_a$ can
query among its first $V_a$ queries, for each $a<j$. Let
\[
\eta=\theta=\frac{e_j}{100\,j\,2^K},
\]
and let $\tau$ and the \emph{pool program} be given by
Lemma~\ref{lem:poolV} for $\ell=j$, $\eta$ and $z$. For a sufficiently large
\emph{trial rank} $r$ let $(G,V)=(G(r),V(r))$ be the output of the pool
program on input $M=r^8$, and let $(S^{(q)},B^{(q)})_{q\le J}$, where
$J=r^{24}$, be the list of \cite[Lemma 4]{LM} for $r,\theta,\tau$, with its
column indices identified with the positions of $G$ in increasing order.
Each $S^{(q)}\subseteq G$ has $r^2$ elements, and each $B^{(q)}$ has rank
$r$. The list is found from $r,\theta,\tau$ by a finite search. For a pair
$(S,B)$ from the list, write $R_i,T_i,W_i$ for the quantities
$R_{ij},T_{ij},W_{ij}$ computed with $(S,B)$ in place of $(S_j,B_j)$.

\begin{fact}[{\cite[Section 5]{LM}}]\label{fact:candidates}
For all sufficiently large $r$,
\begin{multline}\label{eq:average}
\frac1J\sum_{q\le J}\ \sum_{i\le j}\Bigl(\E_{P_{j-1}}A_i(S^{(q)})
+2^K\kappa\bigl(L_i\cap\{T_i(S^{(q)},B^{(q)})<\infty\}\bigr)\Bigr)\\
\le\ \sqrt2\,j\,r^{-2}+2^Kj(\eta+\theta)\ <\ e_j .
\end{multline}
Call $q$ \emph{suitable} if the inner sum in~\eqref{eq:average} is less than
$e_j$. If $q$ is suitable and $(S_j,B_j)=(S^{(q)},B^{(q)})$, then
\textup{(8)} and \textup{(9)} hold for all $i\le j$, whatever the horizon
and the syndrome $y_j$ are, and
\begin{equation}\label{eq:good}
\sum_{i\le j}\E_{P_{j-1}}A_i(S_j)
+\sum_{i\le j}P_{j-1}\bigl(L_i\cap\{T_{ij}<\infty\}\bigr)<e_j .
\end{equation}
\end{fact}

\subsection*{The horizon}
Fix $g$ as in Theorem~\ref{thm:hyp}. We carry out the construction of
\cite{LM} with two changes at each stage: the horizon is the one produced by
the pool program, and the rank is chosen with regard to $g$.

In \cite{LM} the horizon $V_j$ is found by an unbounded search that uses an
advice number $b_j$ approximating a limiting probability. With wrong advice
the search need not terminate. The next lemma shows that the horizon $V(r)$
of the pool program can be used instead. This removes the advice, and it
makes the reconstruction program below halt on every candidate index, which
is what gives a computable bound on its running time.

Fix a stage $j$ and the data of the previous subsection, a sufficiently
large trial rank $r$, and a suitable $q$. Let $(S,B)=(S^{(q)},B^{(q)})$ and
$V=V(r)$. For $y\in\F_2^r$ let
$P^y=P_{j-1}(\,\cdot\mid B(X{\restriction}S)=y)$ and
\[
\Phi(y)=\sum_{i\le j}\Bigl(\E_{P^y}W_i(V)+P^y\bigl(0<R_i(V)<r\bigr)\Bigr).
\]

\begin{lemma}\label{lem:horizon}
Let $y^*$ be the lexicographically first minimizer of $\Phi$. Then
$(S_j,B_j,y_j)=(S,B,y^*)$, $r_j=r$ and $V_j=V$ satisfy
\textup{(6)--(9)} for every $i\le j$.
\end{lemma}

\begin{proof}
(8) and (9) hold by Fact~\ref{fact:candidates}. For $i\le j$ let $N_i$ be the
set of answer sequences $Z\in H_i$ along which $\pi_i$ does not query all of
$G$ among its first $V$ queries. By~\eqref{eq:high} and the density bound,
\begin{equation}\label{eq:Ni}
P_{j-1}(N_i)\le2^K\kappa(N_i)<2^K\eta=\frac{e_j}{100\,j}.
\end{equation}
We claim that, as events for the rule $\pi_i$,
\begin{equation}\label{eq:incl}
\{V<T_i<\infty\}\ \cup\ \{0<R_i(V)<r\}\ \subseteq\ N_i\cup\bigl(L_i\cap\{T_i<\infty\}\bigr).
\end{equation}
Consider an answer sequence in $H_i$. If $V<T_i$, there is no rank drop
among the first $V$ queries, so $R_i(V)=r\ge1$. If $0<R_i(V)<r$, then
$R_i(V)\ge1$ as well. In both cases some column of $B$ is unread after $V$
queries, so some position of $S\subseteq G$ has not been queried, and the
answer sequence is in $N_i$. Now consider an answer sequence in $L_i$. In the
first case $T_i<\infty$ by assumption. In the second case $R_i(V)<r$ means
that a rank drop has occurred, so $T_i\le V$. This proves~\eqref{eq:incl}.

From \eqref{eq:incl}, \eqref{eq:Ni} and \eqref{eq:good},
\begin{equation}\label{eq:tails}
\sum_{i\le j}P_{j-1}(V<T_i<\infty)<\frac{e_j}{100}+e_j<2e_j,
\qquad
\sum_{i\le j}P_{j-1}\bigl(0<R_i(V)<r\bigr)<2e_j .
\end{equation}

\emph{Proof of \textup{(7)}.} The argument of \cite[p.~11]{LM} applies to any
horizon. Write $C_{\pi_i}(w)$ for the set of inputs on which $\pi_i$ has
transcript $w$. The event $\{V<T_i<\infty\}$ is the disjoint union of the
sets $C_{\pi_i}(w)$ over the transcripts $w$ of length at least $V$ that end
immediately before the first rank drop. On each of them the unread columns
of $B$ have rank $r$, so
$P_{j-1}(C_{\pi_i}(w)\cap\{B(X{\restriction}S)=y\})=2^{-r}P_{j-1}(C_{\pi_i}(w))$
for every $y$. Hence $P^y(V<T_i<\infty)=P_{j-1}(V<T_i<\infty)$ for every
$y$, and (7) follows from the first part of~\eqref{eq:tails}.

\emph{Proof of \textup{(6)}.} The set $S$ is disjoint from the earlier
supports and $B$ has full row rank, so under $P_{j-1}$ the syndrome
$B(X{\restriction}S)$ is uniformly distributed on $\F_2^r$. Therefore the
average of $\Phi$ over all $y\in\F_2^r$ equals
\begin{multline*}
\sum_{i\le j}\Bigl(\E_{P_{j-1}}W_i(V)+P_{j-1}\bigl(0<R_i(V)<r\bigr)\Bigr)\\
\le\sum_{i\le j}\E_{P_{j-1}}A_i(S)+\sum_{i\le j}P_{j-1}\bigl(0<R_i(V)<r\bigr)
<e_j+2e_j ,
\end{multline*}
by $W_i(V)\le A_i(S)$, \eqref{eq:good} and~\eqref{eq:tails}. So
$\Phi(y^*)<3e_j$. Every summand of $\Phi(y^*)$ is non-negative, and
$P_j=P^{y^*}$, which gives (6) for each $i\le j$.
\end{proof}

\subsection*{The reconstruction program and its running time}
As in \cite{LM}, the \emph{header} of stage $j$ consists of the earlier
systems $(S_a,B_a,y_a)$, $a<j$, of programs for $\pi_i$ and $h_i$, $i\le j$,
of the pool program, and of $j,z,\eta,\tau,e_j$. It does not depend on the
trial rank. Let $Q_j$ be the following procedure with input $(r,q)$.
\begin{enumerate}[label=\arabic*.]
\item Run the pool program on input $M=r^8$ to obtain $(G,V)$.
\item Find the list of \cite[Lemma 4]{LM} for $r,\theta,\tau$ by exhaustive
search.
\item If $q>r^{24}$, stop. Otherwise let $(S,B)$ be the $q$th entry of the
list, transported to $G$.
\item Compute $\Phi(y)$ for every $y\in\F_2^r$. These are rational numbers,
obtained by enumerating the transcripts of length $V$ of
$\pi_1,\dots,\pi_j$, as in \cite[p.~11]{LM}.
\item Output $(S,B,y^*)$, where $y^*$ is the lexicographically first
minimizer of $\Phi$.
\end{enumerate}
Steps 1 and 2 terminate for all sufficiently large $r$, and steps 3--5 are
finite computations for every $q\le r^{24}$. Let $p_j(r,q)$ be the program
that consists of the header, the procedure $Q_j$, and $r$ and $q$ in
self-delimiting binary notation. As in \cite{LM} there are constants $C_j$
and $c_0$, with $C_j$ independent of $r$ and $q$, such that
\begin{equation}\label{eq:length}
|p_j(r,q)|\le C_j+c_0\log(r+2)\qquad\text{for all }q\le r^{24}.
\end{equation}

Fix a notion of running time for the universal interpreter, so that it is
decidable whether a program halts within $s$ steps with a given output.

\begin{lemma}\label{lem:time}
There are $r^\circ_j\in\N$ and a total computable function $f_j$ such that
for all $r\ge r^\circ_j$ and all $q\le r^{24}$ the program $p_j(r,q)$ halts
within $f_j(r)$ steps.
\end{lemma}

\begin{proof}
Choose $r^\circ_j$ such that steps 1 and 2 terminate for all
$r\ge r^\circ_j$. For such $r$ let $f_j(r)$ be the largest running time of
$p_j(r,q)$ over $q\le r^{24}$, and let $f_j(r)=0$ for $r<r^\circ_j$. The
number $r^\circ_j$ is a constant, so $f_j$ is computable.
\end{proof}

By Lemma~\ref{lem:horizon}, if $r$ is sufficiently large and $q$ is suitable,
then the output of $p_j(r,q)$ and the horizon $V(r)$ satisfy the estimates of
stage $j$.

\subsection*{Choosing the ranks}
We define the blocks by recursion on $j$. At stage $j$ the earlier stages
determine the header, and with it $r^\circ_j$, $f_j$ and $C_j$. There is
$r^*_j\ge r^\circ_j$ such that every $r\ge r^*_j$ has the following two
properties:
\begin{enumerate}[label=\textup{(\roman*)}]
\item the conclusions of Lemma~\ref{lem:poolV}, Fact~\ref{fact:candidates}
and Lemma~\ref{lem:horizon} hold for the trial rank $r$; in particular a
suitable $q$ exists;
\item $r\ge2\bigl(C_j+c_0\log(r+2)\bigr)+j+1$.
\end{enumerate}
Since $g$ is not dominated by $f_j$, there are infinitely many $r$ with
\begin{enumerate}[label=\textup{(\roman*)},resume]
\item $g(r)\ge f_j(r)$.
\end{enumerate}
Let $r_j$ be the least $r\ge r^*_j$ satisfying (iii), let $q_j$ be the least
suitable $q$ for this rank, and let
\[
p_j=p_j(r_j,q_j),\qquad (S_j,B_j,y_j)=\text{the output of }p_j,\qquad V_j=V(r_j).
\]
The support $S_j$ lies above $z$, so it is disjoint from the earlier supports
and hypothesis~\ref{avoid} of Theorem~\ref{thm:LM6} is maintained. By
Lemma~\ref{lem:horizon}, (6)--(9) hold at stage $j$ for all $i\le j$.
Moreover, by \eqref{eq:length}, (ii), Lemma~\ref{lem:time} and (iii),
\begin{equation}\label{eq:pj}
r_j\ge2|p_j|+j+1\qquad\text{and}\qquad p_j\text{ halts within }g(r_j)\text{ steps.}
\end{equation}

\subsection*{A time-bounded test}
For $m,r\ge1$ let $U^g_{m,r}$ be the union of all $E(S,B,y)$ such that
$(S,B,y)$ is a system of rank exactly $r$ that is output within $g(r)$ steps
by some program $p$ with $2|p|+m+1\le r$. Let
$U^g_m=\bigcup_{r\ge1}U^g_{m,r}$. Then $U^g_m\subseteq U_m$.

\begin{lemma}\label{lem:test}
$(U^g_{m+2})_{m\ge1}$ is a Schnorr test relative to $g$.
\end{lemma}

\begin{proof}
There are fewer than $2^{(r-m+1)/2}$ binary strings $p$ with
$2|p|+m+1\le r$, and none if $r\le m$. Each contributes a clopen set of
measure $2^{-r}$ to $U^g_{m,r}$ or nothing, and which of the two can be
decided from $g(r)$. Thus $U^g_{m,r}$ is a clopen set, an index for it can be
computed from $m$, $r$ and $g(r)$, and
$\lambda(U^g_{m,r})<2^{(1-m-r)/2}$. Hence $U^g_{m+2}$ is $\Sigma^0_1(g)$
uniformly in $m$, and for every $R\ge m+2$
\[
\lambda\Bigl(\bigcup_{r>R}U^g_{m+2,r}\Bigr)
<\sum_{r>R}2^{(-1-m-r)/2}=(2+\sqrt2)\,2^{-(m+R+2)/2}.
\]
With $R=m+2$ this gives
$\lambda(U^g_{m+2})<(2+\sqrt2)\,2^{-m-2}<2^{-m}$, because
$U^g_{m+2,r}=\emptyset$ for $r\le m+2$. In general it shows that
$\lambda(U^g_{m+2})$ is approximated to within $2^{2-(m+R)/2}$ by the measure
of the clopen set $\bigcup_{r\le R}U^g_{m+2,r}$, which is computable from
$g$. So $\lambda(U^g_{m+2})$ is computable relative to $g$, uniformly in $m$.
\end{proof}

\subsection*{Proof of Theorem~\ref{thm:hyp}}
By the choices made at each stage, the hypotheses of Theorem~\ref{thm:LM6}
hold, so $\rho$-almost every sequence is KL-random. By~\eqref{eq:pj}, for
every $j$ and every $m\le j$ we have $E_j\subseteq U^g_{m,r_j}$. Hence
$E_{m+2}\subseteq U^g_{m+2}$ for every $m$, and
\[
\supp\rho=\bigcap_{j\ge1}E_j\subseteq\bigcap_{m\ge1}U^g_{m+2}.
\]
So every sequence in $\supp\rho$ fails the Schnorr test of
Lemma~\ref{lem:test}. For the second part of the theorem, let $A$ have
hyperimmune degree and let $g\le_TA$ be a function that is not dominated by
any computable function. A Schnorr test relative to $g$ is a Schnorr test
relative to $A$. Since $\rho(\supp\rho)=1$, some KL-random sequence lies in
$\supp\rho$, and it is not Schnorr random relative to $A$. It is then not
computably random relative to $A$, because $\CR^A\subseteq\SR^A$. \qed

\section{The formalization}\label{sec:lean}

The results on base invariance and on the complexity of the separation,
and the hyperimmune-freeness theorem behind Theorem~\ref{thm:lowklr}, are
formalized in Lean~4 with Mathlib in the repository \cite{klr}, which also
contains the formalization of \cite{LM}. Both bases
use total scan rules and uniformly computable real capital, as in
Section~\ref{sec:bet}. The formal statements are about digit streams: with
$\mathtt{pairDigits}\,X=W_X$, the main theorem reads
\begin{quote}\small\ttfamily
theorem exists\_baseFour\_random\_not\_binary\_random :\\
\hspace*{1.5em}$\exists$ X : Cantor,\\
\hspace*{3em}BaseFour.KolmogorovLovelandRandom (BaseFour.pairDigits X) $\wedge$\\
\hspace*{3em}$\neg$KolmogorovLovelandRandom X
\end{quote}
Table~\ref{tab:lean} gives the correspondence with the statements of this
paper; all names are in the namespace \texttt{Klr}.

\begin{table}[ht]
\renewcommand{\arraystretch}{1.3}
\begin{tabular}{@{}>{\raggedright\arraybackslash}p{0.26\textwidth}>{\raggedright\arraybackslash\footnotesize}p{0.70\textwidth}@{}}
\hline
Statement & {\normalsize Lean declaration}\\
\hline
Theorem~\ref{thm:main} & \path{exists_baseFour_random_not_binary_random}\\
Theorem~\ref{thm:eff}(a) & \path{exists_fourthJump_computable_separating_measure}\\
Theorem~\ref{thm:eff}(b) & \path{exists_kolmogorovLovelandRandom_not_martinLofRandom_below_fifthJump}\\
Corollary~\ref{cor:strict} & \path{binary_random_ssubset_baseFour_random}\\
Proposition~\ref{prop:easy}, $b=k=2$ & \path{BaseFour.of_binary_kolmogorovLovelandRandom}\\
Lemma~\ref{lem:sim} & \path{BaseFour.kolmogorovLovelandRandom_pairDigits_iff}\\
Theorem~\ref{thm:split}, even/odd case & \path{kolmogorovLovelandRandom_even_or_odd}\\
Proposition~\ref{prop:stage} & \path{ConfinedStageConstruction.exists_sequence}\\
Lemma~\ref{lem:even} & \path{ParityPools.exists_fixed_partitions_even_pools}\\
Lemmas~\ref{lem:even} and~\ref{lem:transfer} & \path{ParityPools.exists_fixed_partitions_parity_pools}\\
Lemma~\ref{lem:tests} & \path{ParityDescription.not_martinLofRandom_half_of_blocks}\\
(i) and (ii) for one $X$ & \path{ParityConstruction.exists_pairRespectingRandom_halves_not_martinLofRandom}\\
Theorem~\ref{thm:hyp}, computably random form & \path{exists_kolmogorovLovelandRandom_not_computablyRandomRelative}\\
$\Low(\KLR,\CR)$ is hyperimmune-free & \path{hyperimmuneFree_of_klr_implies_computablyRandomRelative}\\
$\Low(\KLR)$ is hyperimmune-free & \path{lowForKLR_implies_hyperimmuneFree}\\
\eqref{eq:pj} for one $X$ & \path{exists_kolmogorovLovelandRandom_timed_descriptions}\\
Theorem~\ref{thm:LM6} & \path{ae_kolmogorovLovelandRandom_of_stage_estimates}\\
Lemma~\ref{lem:poolV} & \path{TimedPool.exists_fixed_partitions}\\
Lemma~\ref{lem:horizon} & \path{TimedStageHorizon.partial_and_tail_lt}, \path{TimedStage.result_of_data}\\
Lemma~\ref{lem:time} & \path{TimedStage.all_candidates_halt}, \path{TimedProgram.exists_runtime_bound}\\
\hline
\end{tabular}
\medskip
\caption{Statements of this paper and their formal counterparts in \cite{klr}.}\label{tab:lean}
\end{table}

The formal development differs from the text in four respects. The
even/odd case of the splitting theorem is proved in the repository and is
not taken from \cite{MMNRS} as a hypothesis. The tests of
Lemma~\ref{lem:tests} check, by a finite truth table, that the event output
by a program depends only on coordinates of the given parity, instead of
inspecting the support of a presented system. The formal main theorem is
existential; the almost-everywhere statement of Theorem~\ref{thm:main} is
established inside its proof for the block law $\rho'$. Proposition~\ref{prop:easy}
is formalized for $b=k=2$ only.

For Theorem~\ref{thm:eff}, $\jump{n}$ is formalized as the $n$-fold
relativized halting problem over the empty oracle, and relative computability
is Mathlib's. Part (b) reads
\begin{quote}\footnotesize\ttfamily
theorem exists\_kolmogorovLovelandRandom\_not\_martinLofRandom\_below\_fifthJump :\\
\hspace*{1.5em}$\exists$ X : Cantor, KolmogorovLovelandRandom X $\wedge$ $\neg$MartinLofRandom X $\wedge$\\
\hspace*{3em}Reduces X (OracleJump.zeroJump 5)
\end{quote}
In part (a) the formal statement provides a function from finite strings to
rationals, computable in $\jump4$, whose values are the cylinder
probabilities. The effective version of the base-$4$ construction of
Theorem~\ref{thm:main} is not part of the formalization.

\subsection*{Lowness}
The repository formalizes the part of Theorem~\ref{thm:hyp} on which
Theorem~\ref{thm:lowklr} rests:
\begin{quote}\footnotesize\ttfamily
theorem exists\_kolmogorovLovelandRandom\_not\_computablyRandomRelative\\
\hspace*{3em}\char`\{A : Oracle\char`\}\ (hA : $\neg$HyperimmuneFree A) :\\
\hspace*{1.5em}$\exists$ X : Cantor, KolmogorovLovelandRandom X $\wedge$
$\neg$ComputablyRandomRelative A X
\end{quote}
Here \texttt{HyperimmuneFree A} says that every total function computable
from $A$ is eventually bounded by a computable function, and
\texttt{ComputablyRandomRelative A X} says that no martingale whose capital is
computable from $A$ succeeds on $X$. Two consequences are stated separately:
every set in $\Low(\KLR,\CR)$ is hyperimmune-free, and every set in
$\Low(\KLR)$ is hyperimmune-free (Table~\ref{tab:lean}). In the formal
definition of $\Low(\KLR)$ the strategies relative to $A$ are total, with
scan rule and capital computable from $A$, as in Section~\ref{sec:bet}. The
martingale that the proof constructs is total, rational-valued and computable
from $A$, so the theorem does not depend on which of the usual conventions
for relativized strategies is chosen.

The formal proof follows Section~\ref{sec:low} up to \eqref{eq:pj}.
Table~\ref{tab:lean} lists the counterparts of Theorem~\ref{thm:LM6} and of
Lemmas~\ref{lem:poolV}, \ref{lem:horizon} and~\ref{lem:time}. The choice of
the ranks is \path{TimedStageConstruction.exists_fixed_partitions_steps}, and
the recursion is \path{TimedStageRecursion.exists_sequence}. The formal proof
differs from the text in three respects. In addition to
$r_j\ge2|p_j|+j+1$ it requires $r_j\ge4|p_j|$, which again holds for all
sufficiently large ranks. As in the formalization of \cite{LM}, the quantity
that is minimized over the syndromes counts all wagers on the support of the
block among the first $V$ queries; it lies between $W_i(V)$ and $A_i(S)$, so
the argument of Lemma~\ref{lem:horizon} applies to it. And the Schnorr test
of Lemma~\ref{lem:test} is replaced by a martingale. For a system $(S,B,y)$
the function $v\mapsto2^{|v|}\lambda(E(S,B,y)\cap[v])$ is a rational
martingale with initial capital $2^{-\rank B}$ that reaches $1$ along every
solution. The systems that are output within $g(r)$ steps by programs of
length at most $r/4$ are grouped by rank, and the groups are combined into a
single martingale that is computable from $g$ and succeeds on every sequence
in $\supp\rho$.

Three things that this paper uses are not formalized. The first is
Lemma~\ref{lem:test}, an elementary measure computation. With it the Schnorr
form of Theorem~\ref{thm:hyp} follows from the formal theorem
\path{exists_kolmogorovLovelandRandom_timed_descriptions}, which provides a
KL-random sequence in $\bigcap_mU^g_m$. The second is the classical theorems
of \cite{N05,Ch,MM} that complete the proof of Theorem~\ref{thm:lowklr}, so
that the equality in that theorem is not itself a formal statement. The third
is the results of \cite{KNS,TZ} behind Corollary~\ref{cor:klrsr} and
Proposition~\ref{prop:mlrklr}.

\section{Remarks and questions}

\begin{question}\label{q:four}
Is there a sequence that is KL-random, not ML-random, and computable in
$\jump4$? More generally, what is the least $n$ such that some KL-random
sequence that is not ML-random is computable in $\jump{n}$? In particular,
is there a $\Delta^0_2$ such sequence?
\end{question}

Theorem~\ref{thm:eff}(b) gives $n\le5$. The last step of its proof spends
one jump on selecting a point from a measure that is computable in
$\jump4$; an answer with $n=4$ would have to avoid this step or to lower
the complexity of the measure.

\begin{remark}[highness]
Every KL-random sequence that is not ML-random is high, that is,
$X'\ge_T\jump2$. Indeed, KL-random sequences are Schnorr random, and Nies,
Stephan and Terwijn \cite{NST} showed that a Schnorr random sequence that is
not high is ML-random. They also showed that every high degree contains a
computably random sequence that is not ML-random. We do not know which high
degrees contain a KL-random sequence that is not ML-random.
\end{remark}

\begin{question}\label{q:measure}
What is the least $n$ such that there is a probability measure with
cylinder probabilities computable in $\jump{n}$ for which almost every
sequence is KL-random and not ML-random?
\end{question}

By Remark~\ref{rem:notcomp} and Theorem~\ref{thm:eff}(a), $1\le n\le4$.

\begin{remark}[Kurtz randomness]
Let $\WR$ be the class of Kurtz random (weakly random) sequences. Greenberg
and Miller \cite{GM} showed that $\Low(\MLR,\WR)$ consists of the sets that
are not of DNC degree, and that $\Low(\CR,\WR)=\Low(\SR,\WR)$ consists of the
sets that are neither high nor of DNC degree. Hence
\[
\{A:A\text{ not high, not of DNC degree}\}
\ \subseteq\ \Low(\KLR,\WR)\ \subseteq\ \{A:A\text{ not of DNC degree}\}.
\]
For $\CR$ in place of $\KLR$ no high set qualifies, because every high set
computes a computably random sequence \cite{NST}, which is then not Kurtz
random relative to that set. This argument is not available for $\KLR$: by
\cite{MMNRS} both halves of a $\Delta^0_2$ KL-random sequence are ML-random,
so a $\Delta^0_2$ set that is not of DNC degree computes no KL-random
sequence. The method of Section~\ref{sec:low} does not apply directly either.
A test for Kurtz randomness relative to $A$ would require an $A$-computable
bound on the ranks $r_j$, and these depend on the non-effective data of the
construction.
\end{remark}

\begin{question}\label{q:kurtz}
Is there a high set $A$, not of DNC degree, such that every KL-random
sequence is Kurtz random relative to $A$?
\end{question}

\begin{remark}[effective bounds in base $4$]
We expect the bounds of Theorem~\ref{thm:eff} to hold for the measure
$\rho'$ and for the witness of Theorem~\ref{thm:main} as well, by confining
the reservoirs of Section~\ref{sec:eff} to positions of the required parity.
We have not carried this out.
\end{remark}

\begin{remark}[bases $2$ and $2^k$]
For $k\ge2$ the same proof gives a real that is KL-random in base $2^k$ and
not in base $2$; this case is not part of the formalization. Use the binary
strategies that read each block $\{kn,\dots,kn+k-1\}$ consecutively, put the
blocks of even stages on positions $\equiv0$ and those of odd stages on
positions $\equiv1$ modulo $k$, and apply Theorem~\ref{thm:split} with
$Z=k\N$. The subsequence selected by $\N\setminus Z$ is not ML-random because
its computable subsequence on the positions $\equiv1$ modulo $k$ is not.
\end{remark}

\begin{remark}[bases $4$ and $8$]
We expect bases $2^a$ and $2^b$ to be incomparable when neither exponent
divides the other. For $4$ and $8$, group the bits in sixes
$b_0\cdots b_5$. The bits $b_1,b_2$ lie in one base-$8$ digit but in two
base-$4$ digits, and $b_2,b_3$ lie in one base-$4$ digit but in two base-$8$
digits. Constraining these two bit-coordinates at alternating stages should
give KL-randomness in the base whose digits keep the two coordinates
together. Refuting KL-randomness in the other base requires
Theorem~\ref{thm:split} for sequences over a four- or eight-letter alphabet.
\end{remark}

\begin{question}
Is it true that KL-randomness in base $b$ implies KL-randomness in base $b'$
if and only if $b'$ is a power of $b$?
\end{question}

\begin{remark}[multiplicatively independent bases]
For bases such as $2$ and $3$ the method gives nothing. Neither kind of
strategy simulates the other, and the answer to a ternary query is not a
function of boundedly many bits. Under a measure defined by parity
constraints on bits, the answer to a ternary query is therefore no longer
either conditionally fair or forced, which is the dichotomy on which
\cite[Sections 3 and 6]{LM} rest. One needs a quantitative substitute, for
instance a bound on the correlation between a parity of sparse, scattered
bits and the events determined by a ternary transcript. We conjecture that
invariance fails here as well, in both directions.
\end{remark}

\section*{Acknowledgments}
The proofs of the main result were sketched by Claude Fable 5.2 and ChatGPT 6 Pro, and formalized
in Lean by ChatGPT; see \cite{klr}.
This work
was supported by a grant from the Simons Foundation (\#4508914 to Bj{\o}rn Kjos-Hanssen).

\end{document}